\documentclass[11pt, a4paper]{amsart}
\usepackage{a4wide}
\usepackage{amssymb}
\usepackage{amsthm}
\usepackage{amsfonts}
\newtheorem{lemma}{\bf Lemma}[section]
\newtheorem{theorem}[lemma]{\bf Theorem}
\newtheorem{corollary}[lemma]{\bf Corollary}
\newtheorem{proposition}[lemma]{\bf Proposition}
\newtheorem{example}[lemma]{\bf Example}

\newtheorem{definition}[lemma]{\bf Definition}

\title{Remarks on affine complete distributive lattices}
\author{Dominic van der Zypen}
\begin{document}
\maketitle
\parindent = 0mm
\parskip = 2 mm
\begin{quote} {\footnotesize {\bf Abstract.} 
We characterise the Priestley spaces corresponding to affine complete
bounded distributive lattices. Moreover we prove that the class of
affine complete bounded distributive lattices is closed under products
and free products. We show that every (not necessarily bounded) distributive
lattice can be embedded in an affine complete one and that $\mathbb{Q} \cap
[0,1]$ is initial in the class of affine complete lattices.} 
\footnote{The author is grateful for financial support from the Swiss National Science Foundation. 
 } \end{quote}	
\section{Affine complete lattices}
A $k$-ary function $f$ on a bounded distributive lattice $L$ is
called \emph{compatible} if for any congruence $\theta$ on $L$ and
$(a_i,b_i)\in \theta$, $(i=1,...,k)$ we always have
$(f(a_1,...,a_k),f(b_1,....b_k))\in \theta$. It is easy to see
that the projections $pr_i: L^k \to L$ are compatible. With
induction on polynomial complexity one shows that every polynomial
function is compatible (see \cite{Plos}). A lattice $L$ is called
\emph{affine complete}, if conversely every compatible function on
$L$ is a polynomial.

\medskip
G. Gr\"atzer \cite{Gr} gave an intrinsic characterization of
bounded distributive lattices that are affine complete:

\begin{theorem} {\rm (\cite{Gr})} A bounded distributive lattice is affine
complete if and only if it does not contain a proper interval that
is a Boolean lattice in the induced order.
\end{theorem}

Note that in particular, no finite bounded distributive lattice
$L$ is affine complete: Let $x\in L$  be an element  distinct from
1. Then $x$ has an upper neighbor, ie, there exists $y \in L$ such
that $[x,y]=\{x,y \}$ which is isomorphic to the 2-element Boolean
lattice.

\begin{example} The bounded distributive lattices $[0,1]$ and
$[0,1] \times [0,1]$ are affine complete.
\end{example}
\begin{proof}
First, take any $x<y$ in $[0,1]$. Then the element
$a=\frac{x+y}{2} \in [x,y]$ has no complement $a'$ in $[x,y]$:
Otherwise we would have $a \wedge a' =x$ which would imply $a'=x$,
but then $a \vee a' = a \neq y$. So $[x,y]$ is not Boolean, whence
$[0,1]$ has no proper Boolean interval.

Secondly, let $(x_1,x_2) < (y_1,y_2) \in [0,1]\times [0,1]$. With
a similar argument as before, the element $$(\frac{x_1+y_1}{2},
\frac{x_2+y_2}{2}) \in [(x_1,x_2), (y_1,y_2)]$$ does not have a
complement in $[(x_1,x_2), (y_1,y_2)]$. Thus, $[0,1]\times [0,1]$
has no proper Boolean interval and is therefore affine complete.
\end{proof}

\section{Priestley duality}

In \cite{Pr70}, Priestley proved that the category $\mathcal{D}$
of bounded distributive lattices with $(0,1)$-preserving lattice
homomorphisms and the category $\mathcal{P}$ of compact totally
order-disconnected spaces (henceforth referred to as {\it
Priestley spaces}) with order-preserving continuous maps are
dually equivalent. (A compact {\it totally order-disconnected
space} $(X;\tau ,\leq )$ is a poset $(X;\leq )$ endowed with a
compact topology $\tau$ such that, for $x$, $y\in X$, whenever
$x\not\geq y$, then there exists a clopen decreasing set $U$ such
that $x\in U$ and $y\not\in U$.) The functor $D:{\mathcal{D}}\to
{\mathcal{P}}$ assigns to each object $L$ of ${\mathcal{D}}$ a
Priestley space $(D(L);\tau (L),\subseteq)$, where $D(L)$ is the
set of all prime ideals of $L$ and $\tau (L)$ is a suitably
defined topology (the details of which will not be required here).
The functor $E:{\mathcal{P}}\to {\mathcal{D}}$ assigns to each
Priestley space $X$ the lattice $(E(X);\cup ,\cap ,\emptyset ,X)$, where $E(X)$
is the set of all clopen decreasing sets of $X$.  

Priestley duality therefore provides us with a ``dictionary'' between
the world of bounded distributive lattices and a certain category of
ordered topological spaces. This is interesting in particular because 
free products of lattices are ``translated'' into products of Priestley
spaces. We will use this fact for showing that the class of affine
complete bounded distributive lattices is closed under free products.

\section{Affine complete Priestley spaces}

\def\latt{{\mathcal E}(X)}
\def\N{\mathbb{N}}

The aim of this section is to characterize the Priestley spaces
corresponding to affine complete distributive (0,1)-lattices. Such
spaces will be called {\it affine complete Priestley spaces}. In
other words, a Priestley space $X$ is affine complete iff $\latt$
is affine complete.

The following theorem provides a rather straightforward
translation of the algebraic concept of affine completeness in
order-topological terms.
\begin{theorem}\label{affinetheorem}
Let $X$ be a Priestley space. Then the following statements are equivalent:
\begin{enumerate}
\item $\latt$ is affine complete.
\item If $U\subseteq V$ are clopen down-sets and $U\neq V$, then the subposet $V\setminus U$
of $X$ is not an antichain, i.e. $V \setminus U$ contains a pair of distinct comparable
elements.
\end{enumerate}
\end{theorem}
\begin{proof}
$(1) \implies (2)$. Suppose $V\backslash U$ is an antichain. Let $C\in [U,V]
\subseteq \latt$. Take $C'=U \cup (V\backslash C)$.

\textsc{Claim:} $C'$ is a clopen down-set of $X$.

It is clear that $C'$ is a clopen subset of $X$ since $V \setminus C=
V \cap (X\setminus C)$. Now, let $c \in C'$ and assume $x < c$. Then if
$c \in U$, we are done, since $U$ is a down-set. Assume
$c \in V\setminus U$. Since $V$ is a down-set, we get $x\in V$, and
the fact that $V\setminus U$ is an antichain tells us that $x$ cannot be a member
of $V\setminus U$. Therefore $x \in U \subseteq C'$ which proves that $C'$ is indeed
a (clopen) down-set.

 Moreover, $C'$ is the complement of
$C$ in $[U,V]$, i.e. $C\cap C'=U$ and $C \cup C'=V$. Because $C$ was arbitrary,
we see that $[U,V]$ is a proper Boolean interval of $\latt$, whence
$\latt$ is not affine complete.

\medskip
$(2) \implies (1)$. Suppose $U \subseteq V$ are distinct clopen down-sets.
By assumption, there are elements $x, y \in V\backslash U$ such that $x < y$.
There is a clopen down-set $A$ with $x \in A$ and $y \notin A$. Consider
$B= (A \cap V) \cup U$. So $B \in [U, V]$ and $y \notin B$. Now we show that $B$ has no
complement in $[U,V]$: Take any $C \in [U,V]$ with $C \cup B=V$. Then
$y \in C$, but since $C$ is a down-set, we have $x \in C$, thus
$x \in (B \cap C) \backslash U$ and $B \cap C \neq U$.
So whatever $C$ we pick, $C$ is no complement for $B$, i.e.
$B$ is not complemented, and consequently
$[U,V]$ is not Boolean. It follows that no proper interval
of $\latt$ is Boolean.
\end{proof}

We can formulate the above result in a more concise way:

\begin{corollary}
A Priestley space $X$ is affine complete if and only if each nonempty open
set contains two distinct comparable points.
\end{corollary}
\begin{proof}
It follows directly from  theorem \ref{affinetheorem} that if each nonempty open
set contains two distinct points that are comparable, then $X$ is affine
complete.

\medskip
Conversely, suppose that $U$ is a nonempty open set which is an antichain, then there exist
open down-sets $C_1, C_2$ such that $\emptyset \neq C_1 \cap (X \backslash
C_2 ) \subseteq U$. Then $[C_1\cap C_2, C_1]$ is a proper interval
such that $C_1 \backslash (C_1\cap C_2) = C_1 \cap(X \backslash C_2)$
is an antichain (as a subset of the antichain $U$). Thus
theorem \ref{affinetheorem} implies that $X$ is not affine complete.
\end{proof}

Note that the proof works exactly the same way if each occurrence of ``open''
is replaced by ``clopen'' (basically because each Priestley space is
zero-dimensional). So we can state as well:

\medskip
\def\cal{\mathcal}

A Priestley space $X$ is affine complete if and only if each nonempty clopen
set contains two distinct comparable points.

\section{Products of affine complete lattices}
We prove in this section that arbitrary products of
affine complete lattices are affine complete. We don't need
Priestley duality  to do this. Priestley duals of affine complete
lattices, i.e. affine complete Priestley spaces, will come into
play when we consider coproducts of affine complete lattices.

\begin{theorem} \label{affineproducttheorem} If $(L_i)_{i \in I}$ is a family of (bounded)
affine complete lattices,
then $\Pi_{i\in I}L_i$ is affine complete.
\end{theorem}
\begin{proof}
\def\produ{\Pi_{i\in I}L_i}
We prove the contrapositive of the theorem. Suppose that $\produ$ is not
affine complete. Then it contains a proper interval $[\xi, \eta]$ that is
Boolean. There exists some $k\in K$ such that $\xi(k)<\eta(k)$. We claim that
$$[\xi(k), \eta(k)] \subseteq L_k$$ is a Boolean interval. Set $x=\xi(k),
y=\eta(k)$. Suppose $l \in [x,y]$ and define $\lambda \in \produ$
by
$$\lambda(i)=\begin{cases}
l & \text{ if } i=k \\
\xi(i) & \text { if } i\neq k\end{cases}$$
Because $[\xi, \eta]$ is Boolean, there
exists $\lambda' \in \produ$ such that $\lambda \wedge \lambda' = \xi$ and
$\lambda \vee \lambda'=\eta$. Thus it is easy to see that $l':=\lambda'(k)$
is the complement of $l \in [x,y]$. Therefore, $[x,y]$ is a proper Boolean
interval of $L_k$ and whence $L_k$ is not affine complete.
\end{proof}
\begin{example} Theorem \ref{affineproducttheorem} implies that $[0,1]^{\N}$ is affine
complete.
\end{example}

\section{Free products of affine complete lattices}
Now we turn our attention to free products of
affine complete bounded distributive lattices; we prove they
are complete. A convenient way to
obtain this result is to dualise the problem into the category of
Priestley spaces. Free products (that is, coproducts) in ${\bf D}$
correspond to products in ${\bf P}$ and vice versa; this is stated
in the following proposition in a more general way.

\begin{proposition}\cite{MacLane}
Let $\mathcal A$ and $\mathcal B$ be categories, and assume that ${\bf
F}:{\mathcal A} \to {\mathcal B}$ and ${\bf G}:{\mathcal B} \to {\mathcal A}$ are
contravariant functors that form a dual equivalence. Then:
\begin{enumerate}
\item If $A$ is a product of a family of objects $(A_i)_{i\in I}$
of $\mathcal{A}$, then ${\bf F}(A)$ is a coproduct of $({\bf
F}(A_i))_{i\in I}$.
\item If $A$ is a coproduct of a family of
objects $(A_i)_{i\in I}$ of $\mathcal{A}$, then ${\bf F}(A)$ is a
product of $({\bf F}(A_i))_{i\in I}$.
\end{enumerate}
\end{proposition}

Moreover we have shown that affine complete lattices correspond
to affine complete spaces under the Priestley duality.

\begin{theorem} If $(X_i)_{i \in I}$ is a family of
affine complete Priestley spaces,
then $\Pi_{i\in I}X_i$ is affine complete.
\end{theorem}
\begin{proof}
\def\produ{\Pi_{i\in I}X_i}
Suppose that $X_i$ is affine complete for every $i\in I$. It suffices to
show that every nonempty subset $V$ of $\produ$ of the form
$$V=\pi^{-1}_{i_1}(U_1) \cup ... \cup \pi^{-1}_{i_r}(U_r)$$ contains
two distinct comparable elements (where $U_{k} \subseteq X_{i_{k}}$ open,
nonempty). Take $U_1$. It contains elements $a<b$, because $X_{i_1}$ is
affine complete. Now pick $\xi \in V$. Define $\xi_1, \xi_2 \in V$ by
$$\xi_1(i)=\begin{cases}
\xi(i) & \text{ if } i\neq i_1 \\
a & \text { if } i=i_1\end{cases}$$
and
$$\xi_2(i)=\begin{cases}
\xi(i) & \text{ if } i\neq i_1 \\
b & \text { if } i=i_1\end{cases}$$
Clearly, $\xi_1, \xi_2$ are distinct comparable elements of $V$.
\end{proof}
Applying the Priestley duality now yields:

\begin{corollary} The class of (bounded) affine complete lattices is closed
under free products.
\end{corollary}

\section{Embedding lattices in  affine complete lattices \label{affineembed}}

\def\L{L_{01}}

First we will stay away from affine completeness in the worst
possible way: we will embed each $L$ into a powerset of some set,
which, being Boolean, is as affine incomplete as it gets.
\begin{lemma}\label{lemma2point8}
Let $L$ be a distributive lattice ($L$ need not be bounded).
There is a set $X$ and
a lattice embedding $$j:L \hookrightarrow {\cal P}(X)$$
where ${\cal P}(X)$ is the powerset of the set $X$.
\end{lemma}
\begin{proof}
First, endow $L$ with a smallest element and a greatest element. Call this new
bounded distributive lattice $\L$. By Priestley duality, there is a Priestley
space $(X,\tau, \leq)$ such that the lattice $\latt$ of clopen down-sets
is isomorphic to $\L$. Since $\latt$ is a sublattice of ${\cal P}(X)$, we are
done.
\end{proof}
Next, we will embed that powerset in
an affine complete lattice.

\begin{lemma}\label{lemma2point9}
Let $X$ be a set and let $Q=\{q \in \mathbb{Q}; 0 \leq q \leq
1\}$. Then there is a lattice embedding
$$j: {\cal P}(X) \hookrightarrow Q^X.$$
Moreover, $Q$ is affine complete.
\end{lemma}
\begin{proof}
Set $j: S \mapsto \chi_S \in Q^X$ for every $S \subseteq X$, where
$\chi_S$ is defined by $$\chi_S(x)=
\begin{cases}
1 & \text{ if } x \in S\\
0 & \text{ if } x \notin S
\end{cases}  $$ It is easy to see that $j$ is a lattice embedding.
Next, we claim that $Q$ is affine complete. Take any $x<y$ in $Q$.
Then the element $a=\frac{x+y}{2} \in [x,y]$ has no complement
$a'$ in $[x,y]$: Otherwise we would have $a \wedge a' =x$ which
would imply $a'=x$, but then $a \vee a' = a \neq y$. So $[x,y]$ is
not Boolean, whence $Q$ has no proper Boolean interval. Therefore,
$Q$ is affine complete.

Moreover, by \ref{affineproducttheorem}, $Q^X$ is affine complete
which concludes the proof.
\end{proof}
Lemmas \ref{lemma2point8} and \ref{lemma2point9} now imply:
\begin{corollary}
Every distributive lattice (not necessarily bounded) can be embedded in a bounded affine complete lattice.
\end{corollary}
Admittedly, the construction provided by \ref{lemma2point8} and
\ref{lemma2point9} is highly non-unique and has no minimality
properties.

\def\Q{\mathbb{Q}_{01}}
\section{$\Q$ as initial object in the category of affine
complete lattices\label{q01embed}} The aim of this section is to
show that the lattice $\Q=\mathbb{Q}\cap [0,1]$ can be embedded
into each affine complete lattice, which amounts to saying that
$\Q$ is an initial object of the category of affine complete
lattices (with (0,1)-homomorphisms, i.e. a full subcategory of the
category bounded distributive lattices). The key will be the
notion of a dense chain.

\begin{definition}
A chain $(X,\leq)$ is called dense if for all $x<y \in X$ there is $z\in X$
with $x<z<y$.
\end{definition}
The first tool we need here is a well known result of model
theory. It states that the theory of dense linear orders is
complete and has $(\mathbb{Q},\leq)$ as prime model. We will state
this result in a more primitive way and prove it.
\begin{proposition}\label{prop4point13}
If $(X,\leq)$ is a bounded dense chain, there is a $(0,1)$-embedding
$$\varphi: \Q \hookrightarrow X.$$
\end{proposition}
\begin{proof}
Let $a:\omega \to \Q\backslash\{0,1\}$ be a bijection. We will
write $a_k$ instead of $a(k)$ to simplify  notation and will
inductively build a subset
$$f \subseteq ( \Q\backslash\{0,1\}) \times (X\backslash\{0_X,1_X\})$$
that's an injective function from $\Q\backslash\{0,1\}$ to
$X\backslash\{0_X,1_X\}$ which is even order-preserving.

\underline{$n=0$:} Choose $b_0\in X\backslash\{0_X,1_X\}$ and set
$f_0:=\{(a_0,b_0)\}$.

\underline{$n\to n+1$:} Assume that $f_n$ has been defined in a
way that for all $k,l \in \{0,...,n\}$ the relation $a_k\leq a_l$
implies $f_n(a_k) \leq f_n(a_l)$ and that $f_n$ is an injective
function from $\{a_0,...,a_n\}$ to $X\backslash\{0_X,1_X\}$. Now
consider the element $a_{n+1} \in \Q \backslash \{0,1\}$.

{\em Case 1:} $a_{n+1}\geq a_i$ for all $i\in \{0,...,n\}$. Then, since
$X$ is dense, there is $b_{n+1} \in X$ such that $1_X>b_{n+1}\geq f_n(a_i)$
for all $i\in \{0,...,n\}$. So, $$f_{n+1}:=f_n\cup \{(a_{n+1},b_{n+1})\}$$
is an injective order-preserving function  that continues  $f_n$.

{\em Case 2:} $a_{n+1}\leq a_i$ for all $i\in \{0,...,n\}$. Proceed similarly
as in Case 1.

{\em Case 3:} There are $k,l\in \{0,...,n\}$ such that $a_k<a_{n+1}<a_l$.
We may assume that there is no $k'\in\{0,...,n\}$ with $a_k<a_{k'}<a_{n+1}$
and likewise that there is no $l'\in \{0,...,n\}$ with $a_{n+1}<a_{l'}<a_l$.
Consider $b_k=f_n(a_k)$ and $b_l=f_n(a_l)$. Since $f_n$ is order-preserving
and injective by assumption, we get $b_k<b_l$. Because $X$ is dense, there
is an element $b_{n+1}$ such that $b_k<b_{n+1}<b_l$. Then
$$f_{n+1}:=f_n\cup \{(a_{n+1},b_{n+1})\}$$
is easily seen to be an injective order-preserving map that
continues $f_n$.

Now, it is easy to see that $$f:=\bigcup_{n\in \omega}f_n$$ is an
injective order-preserving function from $\Q\backslash\{0,1\}$ to
$(X\backslash\{0_X,1_X\}$ which is even order-preserving. So
$$\varphi:=f \cup \{(0,0_X),(1,1_X)\}$$ is an order embedding from
$\Q$ to $X$.
\end{proof}
\begin{proposition} \label{proppie}
Let $L$ be a bounded affine complete distributive lattice. Then
\begin{itemize}
\item[a)] There is a maximal chain $C\subseteq L$, i.e., a chain
that is not properly contained in another chain in $L$. \item[b)]
If $C$ is a maximal chain of $L$ then $C$ is dense.
\end{itemize}
\end{proposition}
\begin{proof}
$a)$ is a standard application of Zorn's Lemma: \def\K{{\cal K}} If $\K$ is
a set of chains of $L$ such that for any $C_1, C_2 \in \K$ we either have
$C_1 \subseteq C_2$ or $C_1\supseteq C_2$, then $\bigcup \K$ is easily checked
to be a chain in $L$: Let $x,y \in \bigcup \K$, then there are members
$C, D$ containing $x,y$ respectively; now since $\K$ is a chain with respect
to $\subseteq$, at least one of the statements  $x,y \in C$ or $x,y \in D$
holds. Since $C,D$ are chains in $L$, either statement leads us
to  $x\leq_L y$ or $x\geq_L y$. So $\K$ is bounded in the poset of all chains
of $L$, thus Zorn's Lemma  implies that there is a maximal chain.

As for $b)$, assume that $C$ is a maximal chain such that $x<y \in
C$ but there is no $z\in C$ with $x<z<y$. Now if there were no $z$
in the whole lattice $L$ such that $x<z<y$, then $[x,y]=\{x,y\}$
is a proper Boolean interval of $L$ which implies that $L$ is not
affine complete,  leading to a contradiction. Thus there is such a
$z$, whence $C\cup\{z\}$ is a chain of $L$ that properly contains
$C$, contradicting the maximality of $C$.
\end{proof}
Now the propositions \ref{prop4point13} and \ref{proppie} directly
imply the following.
\begin{theorem}
If $L$ is an affine complete lattice, then there exists a (0,1)-embedding
$\varphi: \Q \hookrightarrow L.$ \label{q01theorem}
\end{theorem}
\begin{proof} Pick any maximal chain $C$ in $L$. Note that by maximality of $C$
we have $0,1 \in C$ since $C \cup \{0,1\}$ is a chain. So the
inclusion map $\iota: C\hookrightarrow L$ is a $(0,1)$-embedding
as well as the embedding from $\Q$ to $C$ provided by proposition
\ref{proppie}. Composing these two, we get a (0,1)-embedding from
$\Q$ to $L$.
\end{proof}
\section{Open questions}
In chapter \ref{affineembed} we showed that ever bounded
distributive lattice can be extended to an affine complete
lattice. This was achieved by making use of $\Q$ which happens to
be embeddable in any affine complete lattice, ie, the ``smallest''
affine complete lattice. Now the question is: Is the construction
carried out in chapter \ref{affineembed} in some way canonical?
For an arbitrary lattice $L$, does its 'affine hull' have any
interesting universal properties?

\parskip = 0mm

\bigskip

D.~van der Zypen\\
Allianz Suisse Insurance Company\\
Bleicherweg 19\\
CH-8022 Zurich, Switzerland\\
{\tt dominic.zypen@gmail.com}


{\footnotesize
\begin{thebibliography}{99}

\bibitem{LO} B.A.Davey and H.A.Priestley,  {\bf Lattices and Order}, 
Cambridge University Press, 1990.

\bibitem{Gr} G. Gr\"atzer, {\it Boolean functions on distributive lattices},
Universal Algebra and Applications, vol. {\bf 9}, Banach Center Publications,
Warsaw, 97-104.
	
\bibitem{MacLane} S.MacLane, {\bf Categories for the working
mathematician}, Springer Verlag; 2nd edition (1998).

\bibitem{Plos}  M.Plo\v{s}\v{c}ica, {\it Affine Complete Distributive Lattices
}, Order {\bf 11} (1994), 385-390.

\bibitem{Pr70}  H.A.Priestley, {\it Representation of distributive
lattices by means of ordered Stone spaces}, Bull. London Math.
Soc. 2 (1970), 186--190.

\bibitem{H1} H.A.Priestley,  {\it Ordered topological spaces and the representation of distributive lattices}
Proc. London Math. Soc. (3)  {\bf 24}  (1972), 507--530.

\bibitem{vdZ} D.van der Zypen, {\it Aspects of Priestley Duality}, PhD thesis,
	University of Bern, 2004.
\end{thebibliography}
}
\end{document}